\documentclass[a4paper]{article}
\usepackage[english]{babel}
\usepackage[ansinew]{inputenc}
\usepackage{amsmath, amssymb, enumerate, mathtools}
\usepackage[colorlinks=true,linkcolor=black,citecolor=black,filecolor=black]{hyperref}
\usepackage[margin=1.1in]{geometry}

\binoppenalty=\maxdimen 
\relpenalty=\maxdimen
\DeclarePairedDelimiter\floor{\lfloor}{\rfloor}

\usepackage[thmmarks]{ntheorem}
\newtheorem{thm}{Theorem}[section]
\newtheorem{prop}[thm]{Proposition}
\newtheorem{lemma}[thm]{Lemma}
\newtheorem{cor}[thm]{Corollary}
\theorembodyfont{\normalfont}
\newtheorem{defn}[thm]{Definition}
\newtheorem{rmk}[thm]{Remark}
\newtheorem{exa}[thm]{Example}
\newtheorem{con}[thm]{Conjecture}
\theoremstyle{nonumberplain}
\theoremsymbol{\ensuremath{\square}}
\newtheorem{prf}{Proof}

\usepackage[backend=bibtex, style=numeric]{biblatex}
\usepackage{filecontents}

\bibliography{shortbib.bib}

\newcommand{\Cb}{\mathbb{C}}
\newcommand{\Hb}{\mathbb{H}}
\newcommand{\Zb}{\mathbb{Z}}
\newcommand{\Qb}{\mathbb{Q}}
\newcommand{\Nb}{\mathbb{N}}
\newcommand{\slz}{\text{SL}_2(\mathbb{Z})}
\newcommand{\mc}{\mathrm{M}}
\newcommand{\mct}{\widetilde{\mathrm{M}}}
\newcommand{\ft}{\widetilde{f}}
\newcommand{\Dc}{\mathcal{D}}
\newcommand{\Fc}{\mathcal{F}}

\title{On a conjecture of Kaneko and Koike}
\author{Andreas Mono}
\date{} 

\begin{document}


\maketitle

\begin{abstract}
In \cite{kk06}, Kaneko and Koike defined extremal quasimodular forms and proved their existence in depth $1$ and $2$. After normalizing and restricting to the case of depth at most $4$, they conjectured a certain bound on the Fourier coefficients of such forms. More precisely, the prime factors of the denominators of the coefficients are requested to be smaller than the weight. Recently, Pellarin proved this conjecture in the case of depth $1$ and weight divisible by $6$. In this paper, we complete the picture in depth $1$. First, we show that his result implies the same result in the case of weight $6k+4$ for every integer $k \geq 0$ directly. Secondly, we adapt the strategy of his proof in the case of weight $w = 6k$ to the last case of weight $w = 6k+2$. Finally, we provide all computational details to both his and our intermediate results, since those details are essential to his proof, but were omitted during his exposition in \cite[6-10]{pellarin}. Parallel and independent from this work, Peter Grabner proved the aforementioned conjecture in full generality, see \cite[25-32]{grabner}.
\end{abstract}

\bigskip

\textbf{Keywords:}
\textit{Normalized extremal quasimodular forms, Fourier coefficients, bound prime factors.}

\bigskip

\section{Preliminaries and statement of results}

\hspace{\parindent} Throughout this paper we work on the upper half plane $\Hb$ with the group $\Gamma \coloneqq \slz$ acting on it via fractional linear transformations. Recall the definition of the normalized holomorphic Eisenstein series
\begin{align*}
E_2(z) \coloneqq 1 - 24\sum_{n \geq 1} \sigma_1(n)q^n \qquad E_4(z) \coloneqq 1 + 240\sum_{n \geq 1}\sigma_3(n)q^n \qquad E_6(z) \coloneqq 1 - 504\sum_{n \geq 1}\sigma_5(n)q^n
\end{align*}
and $z \in \Hb $, $q \coloneqq \text{e}^{2\pi i z}$, $\sigma_j(n) \coloneqq \sum_{d \mid n} d^j$.

Quasimodular forms extend the theory of classical modular forms naturally if one wishes to preserve holomorphicity. To this end, we relax the transformation law imposed on classical modular forms.
A first motivation towards such a relaxation arises from the fact that the normalized derivative of a classical modular form $f$, given by
\begin{align*}
D(f) \coloneqq \frac{1}{2 \pi i}\frac{df}{dz} = q\frac{df}{dq},
\end{align*}
fails to preserve classical modularity. To see this, we compute
\begin{align*}
D(f)\left(\frac{az+b}{cz+d}\right) = (cz+d)^{w+2}\left(D(f)(z) + \frac{w}{2\pi i}\frac{c}{cz+d}f(z)\right),
\end{align*}
where $w$ denotes the weight of $f$. A second motivation is the transformation law for $E_2$, explicitly
\begin{align*}
E_2\left(\frac{az+b}{cz+d}\right) = (cz+d)^2\left(E_2(z) - \frac{6i}{\pi}\frac{c}{cz+d}\right),
\end{align*}
which is caused by the absence of absolute convergence of the series defining $E_2$. 

Quasimodular forms provide a natural framework for the previous two equations. In other words, $E_2$ is the archetypal example of a quasimodular form that is not a classical modular form. This established, we may consider any positive power of $E_2$ to achieve more generality. The highest occuring power of $E_2$ determines the so called \textit{depth} of a quasimodular form. 

We specify these quantities in the following definition.
\begin{defn} \label{qmfdef}
Let $\ell,w \geq 0$ be integers. A \textit{quasimodular form of weight $w$ and depth $\ell$} is a function $f\colon \Hb \to \Cb$ satisfying one of the two equivalent conditions:
\begin{enumerate}[(i)]
\item The function $f$ is a homogeneous polynomial of degree $\ell$ in the ring $\Cb[E_2,E_4,E_6]$. In other words, $f$ may be written uniquely in the following fashion:
\begin{align*}
f(z) = \sum_{j=0}^{\ell} f_j(z)E_2(z)^j \quad , \quad f_{\ell} \neq 0,
\end{align*}
where $f_j \in \mc_{w-2j}$, especially $f_j$ is an usual modular form of weight $w-2j$.
\item For fixed $z \in \Hb$ and varying $\gamma = \left(\begin{smallmatrix} a & b \\ c & d \end{smallmatrix}\right) \in \Gamma$ the function $f$ has the transformation property
\begin{align*}
\left(f\vert_w\gamma\right)(z) \coloneqq f(\gamma z) (cz+d)^{-w} = P_f\left(\frac{c}{cz+d}\right)
\end{align*}
for some ordinary\footnote{meaning 'not necessarily homogeneous'} polynomial $P_f$ of degree $\ell$ and with coefficients $f_j$ as in (i).
\end{enumerate}
We denote the vector space of all quasimodular forms of weight $w$ and depth $\ell$ by $\mct_{\ell,w}$.
\end{defn}

\begin{rmk}
Indeed both definitions are equivalent by the well-known graded direct sum
\begin{align*}
\Cb[E_4, E_6] = \bigoplus_{w \in 2\Nb_0} \mc_w
\end{align*}
and the aforementioned transformation law for $E_2$. In particular it follows that $\mct_{0,w} = \mc_w$ and hence $P_f(0) = f$.
\end{rmk}

In particular, according to Definition \ref{qmfdef} (ii) , we obtain
\begin{align*}
D \colon \mct_{\ell,w} \to \mct_{\ell+1,w+2}
\end{align*}
as desired. This yields a third characterization of $\mct_{\ell,w}$, explicitly given by
\begin{align*}
\mct_{\ell,w} = \begin{cases}
\bigoplus_{j=0}^{\ell} D^j(\mc_{w-2j}) & \text{ if } \ell < \frac{w}{2}, \\
\bigoplus_{j=0}^{\frac{w}{2}-1} D^j(\mc_{w-2j}) \oplus \Cb D^{\frac{w}{2}-1}(E_2) & \text{ if } \ell \geq \frac{w}{2},
\end{cases}
\end{align*}
A proof can be found in \cite[59]{123}.

Two immediate consequences of Definition \ref{qmfdef} (i) are
\begin{align*}
\text{dim}_{\Cb}\left(\mct_{\ell,w}\right) = \sum_{j=0}^{\ell} \text{dim}_{\Cb}(\mc_{w-2j}) \ \text{, where } \text{dim}_{\Cb}(\mc_{k}) = 
\begin{cases} 
0 & \text{ if } k < 0, \\
\floor*{\frac{k}{12}} & \text{ if } k \geq 0 \text{ and } k \equiv 2 \pmod{12}, \\
\floor*{\frac{k}{12}}+1 & \text{ else,}
\end{cases}
\end{align*}
and the chain of inclusions
\begin{align*}
\mc_w = \mct_{0,w}\subseteq \mct_{1,w} \subseteq \ldots \subseteq \mct_{\frac{w}{2}-2,w} \stackrel{\mc_2=\{0\}}{\ = \ } \mct_{\frac{w}{2}-1,w} \subseteq \mct_{\frac{w}{2},w} = \mct_{\frac{w}{2}+1,w} = \ldots \ \ .
\end{align*}
Thus, we stipulate $0 \leq \ell \leq \frac{w}{2}$ and $\ell \neq \frac{w}{2}-1$. 

\begin{rmk}[\protect{\cite[2-6]{pellarin}}]
Restricting ourselves to depth $\ell \leq 4$, the dimension formula reduces to $\text{dim}_{\Cb}(\mct_{\ell,w}) = \text{dim}_{\Cb}(\mc_{(\ell+1)w})$. Recall that there exists a basis $f_1,\ldots,f_{\text{dim}_{\Cb}(\mc_{w})}$ of $q$-expansions for $\mc_w$ for any weight $w$, such that if $a_n(f_j)$ denotes the $n^{\textit{th}}$ coefficient of $f_j$ then $a_n(f_j) = \delta_{n,j}$. (This basis is called the Miller basis for $\mc_{w}$.) Hence, roughly speaking, it is possible to ``diagonalize'' $\mct_{\ell,w}$ in this case. Precisely, assuming $\ell \leq 4$, there exists a basis of $\mct_{\ell,w}$ such that the $q$-expansions of the elements of $\mct_{\ell,w}$ coincide with the canonical diagonal basis $(1, q, q^2, \ldots)$ of $\Cb[[q]]$ on its first $\text{dim}_{\Cb}(\mct_{\ell,w})$ summands.
\end{rmk}

A guiding principle in mathematics is that a theory becomes richer by imposing more structure on it. Applied to our setting, Kaneko and Koike gave the following definition, motivated further by the previous remark.
\begin{defn}
Let $f \in \mct_{\ell,w} \setminus \mct_{\ell-1,w}$ and $m \coloneqq \text{dim}_{\Cb}(\mct_{\ell,w})$. Write $f(z) = \sum_{n \geq 0} a_n q^n$. Then $f$ is called \textit{extremal} if
\begin{align*}
a_0 = \ldots = a_{m-2} = 0 \quad \text{ and } \quad a_{m-1} \neq 0,
\end{align*}
and is called \textit{normalized} if additionally $a_{m-1} = 1$. We abbreviate $\ft_{\ell,w}$ for the normalized extremal quasimodular form of depth $\ell$ and weight $w$.
\end{defn}

We require \textit{Ramanujan's differential system} to construct some examples. It is given by the three equations
\begin{align*}
D(E_2) = \frac{E_2^2-E_4}{12}, \qquad D(E_4) = \frac{E_2E_4-E_6}{3}, \qquad D(E_6) = \frac{E_2E_6-E_4^2}{2}.
\end{align*}

Now, the following examples of normalized extremal quasimodular forms in low weights and low depths become apparent.
\begin{exa}[\protect{\cite[459]{kk06}}] \label{ex1}
We have
\begin{align*}
\ft_{1,2}(z) &= E_2(z) = 1 - 24q - 72q^2 \pm \ldots ,\\
\ft_{1,6}(z) &= \left(\frac{D(E_4)}{240}\right)(z) = \left(\frac{E_2E_4-E_6}{720}\right)(z) = q + 18q^2 + 84q^3 \pm \ldots ,\\
\ft_{1,8}(z) &= \left(-\frac{D(E_6)}{504}\right)(z) = \left(\frac{E_4^2-E_2E_6}{1008}\right)(z) = q + 66q^2 + 732q^3 \pm \ldots ,\\
\ft_{2,4}(z) &= \left(-\frac{D(E_2)}{24}\right)(z) = \left(\frac{E_2^2-E_4}{288}\right)(z) = q + 6q^2 + 12q^3 \pm \ldots ,\\
\ft_{2,8}(z) &= \left(\frac{5E_4^2+2E_2E_6-7E_2^2E_4}{362880}\right)(z) = q^2 + 16q^3 + 102q^4 \pm \ldots ,\\
\ft_{3,6}(z) &= \left(\frac{5E_2^3-3E_2E_4-2E_6}{51840}\right)(z) = q^2+8q^3+30q^4 \pm \ldots ,\\
\ft_{4,8}(z) &= \left(\frac{5E_4^2+16E_2E_6+14E_2^2E_4-35E_2^4}{11612160}\right)(z) = q^3+\frac{21}{2}q^4 + 54q^5 \pm \ldots \ .
\end{align*}
\end{exa}

The example $\ft_{4,8}$ demonstrates that there exist normalized extremal quasimodular forms with rational but non-integer coefficients in the case of depth $4$. In view of more numerical experiments, Kaneko and Koike offered the following conjecture.
\begin{con}[\protect{\cite[469]{kk06}}] \label{kkcon}
If the depth is at most $4$, then the Fourier coefficients of any normalized extremal quasimodular form of weight greater than $2$ are always positive. Moreover, no denominator of such coefficients has prime factors greater than the weight.
\end{con}

Invoking a common partial fraction approach $\frac{1}{\prod_j (x-p_j)} = \sum_j \frac{\xi_j}{x-p_j}$, evaluated at $x = 0$, we investigate the question wether
\begin{align*}
\ft_{\ell,w}(z) \in \Zb\left[\frac{1}{p} \colon p < w\right][[q]]
\end{align*}
depending on given weight $w$ and depth $\ell$. Formally, $\Zb\left[\frac{1}{p}\right]$ is defined as the localization of $\Zb$ at the prime ideal $(p)$, but equivalently we will regard $\Zb\left[\frac{1}{p}\right]$ as ring of polynomials in $p^{-1}$ over $\Zb$.

A new approach to attack this conjecture in depth $1$ was pioneerd by Pellarin in \cite{pellarin}, where he proved the following Theorem.
\begin{thm}[\protect{\cite[Theorem 3.3]{pellarin}}] \label{pellthm}
Let $k \geq 0$ be an integer and $\ft_{1,6k}$ be the normalized extremal quasimodular form of weight $6k$ and depth $1$. Then $\ft_{1,6k} \in \Zb\big[\frac{1}{p} \colon p < 6k\big][[q]]$.
\end{thm}

Consequently, we split the possible weights into residue classes modulo $6$. In other words, we write 
\begin{align*}
w = 6k \quad \text{or} \quad w = 6k+2 \quad \text{or} \quad w = 6k+4.
\end{align*}
Our results verify the analogous claims in the remaining residue classes of the weight and in depth $1$.

\begin{thm} \label{INITHM} 
Let $k \geq 0$ be an integer and $\ft_{1,6k+2}$ be the normalized extremal quasimodular form of weight $6k+2$ and depth $1$. Then $\ft_{1,6k+2} \in \Zb\big[\frac{1}{p} \colon p < 6k+2\big][[q]]$.
\end{thm}

We have a slightly stronger bound than the weight itself in the case of weight $6k+4$.

\begin{thm} \label{CONG4CASE} 
Let $k \geq 0$ be an integer and $\ft_{1,6k+4}$ be the normalized extremal quasimodular form of weight $6k+4$ and depth $1$. Then $\ft_{1,6k+4} \in \Zb\big[\frac{1}{p} \colon p < 6k\big][[q]]$.
\end{thm}

\begin{rmk}
We emphasize once more that Peter Grabner proved Conjecture \ref{kkcon} in full generality first, see \cite[32]{grabner}. However, we would like to present a slightly different approach in the case of depth $1$.
\end{rmk}

The paper is organized as follows. We begin with the proof of Theorem \ref{CONG4CASE} in the upcoming section. Then we move to a short survey of some additional properties of extremal quasimodular forms of depth $1$. The subsequent section \ref{d1sec} is devoted to the proof of Theorem \ref{INITHM} and divided into two parts: The purpose of Subsection \ref{stratsubsec} is to present the main steps of the proof by adapting \cite[6-10]{pellarin}. Subsection \ref{compsubsec} collects all computational details to both Pellarin's and our intermediate steps of the proof. A short outlook is the content of section \ref{outsec}, which concludes the present paper.

\subsection*{Acknowledgements:} 
\hspace{\parindent} The author would like to thank his PhD-advisor Professor Dr.\ Kathrin Bringmann for suggesting Pellarin's paper as well as for her continuous feedback to the work on it. Furthermore, the author would like to thank Joshua Males for many helpful comments on an earlier version of this paper.

\section{Proof of Theorem \ref{CONG4CASE}}
\hspace{\parindent} We prepare the proof by the following result.
\begin{prop} \label{baseprop}
Assume the notation from the previous section.
\begin{enumerate}[(a)]
\item It holds that
\begin{align*}
\text{\normalfont dim}_{\Cb}\left(\mct_{1,6k}\right) = \text{\normalfont dim}_{\Cb}\left(\mct_{1,6k+2}\right) = \text{\normalfont dim}_{\Cb}\left(\mct_{1,6k+4}\right) = k + 1.
\end{align*}
Hence,
\begin{align*}
\ft_{1,w}(z) = q^k + O\left(q^{k+1}\right).
\end{align*}
\item We have that
\begin{align*}
E_4\cdot\mct_{1,6k} = \mct_{1,6k+4}.
\end{align*}
\end{enumerate}
\end{prop}

\begin{prf}
We prove item a). First, observe that both $6k$ and $6k+4$ are never congruent to $2 \pmod{12}$. Thus,
\begin{align*}
\text{dim}_{\Cb}\left(\mct_{1,6k+4}\right) = \text{dim}_{\Cb}\left(\mc_{6k+4}\right) + \text{dim}_{\Cb}\left((\mc_{6k+2}\right) = \text{dim}_{\Cb}\left(\mc_{6k}\right) + \text{dim}_{\Cb}\left(\mc_{6k+2}\right) = \text{dim}_{\Cb}\left(\mct_{1,6k+2}\right).
\end{align*}
If $k = 2j$ is even then $6k+2 \equiv 2 \pmod{12}$, giving
\begin{align*}
\text{dim}_{\Cb}\left(\mct_{1,6k+2}\right) = \floor*{\frac{k}{2}} + 1 + \floor*{\frac{k}{2}+\frac{1}{6}} = 2j+1 = k+1.
\end{align*}
If $k = 2j+1$ is odd, then $6k+2 \not\equiv 2 \pmod{12}$, giving
\begin{align*}
\text{dim}_{\Cb}\left(\mct_{1,6k+2}\right) = \floor*{j+\frac{1}{2}} + 1 + \floor*{j+\frac{1}{2}+\frac{1}{6}} + 1 = 2j+2 = k+1
\end{align*}
as well. Next,
\begin{align*}
\text{dim}_{\Cb}\left(\mct_{1,6k}\right) = \text{dim}_{\Cb}\left(\mc_{6k}\right) + \text{dim}_{\Cb}\left(\mc_{6k-2}\right) = \floor*{\frac{k}{2}} + 1 + \floor*{\frac{k}{2}-\frac{1}{6}} + 1.
\end{align*}
Writing $k=2j+\varepsilon$, $\varepsilon \in \{0,1\}$, we deduce
\begin{align*}
\text{dim}_{\Cb}\left(\mct_{1,6k}\right) = j + 1 + \floor*{j+\frac{\varepsilon}{2}-\frac{1}{6}} + 1 = 
\begin{cases} j + 1 + (j-1) + 1  = 2j+1 = k+1 & \text{ if } \varepsilon = 0, \\
j + 1 + j + 1 = 2j+\varepsilon+1 = k+1 & \text{ if } \varepsilon = 1,
\end{cases}
\end{align*}
and this verifies a). \\
Now, b) is immediate. Supposing $f = f_0 + f_1E_2 \in \mct_{1,6k}$, where $f_0 \in \mc_{6k}$, $f_1 \in \mc_{6k-2}$, we see that $E_4f_0 \in \mc_{6k+4}$, $E_4f_1 \in \mc_{6k+2}$ and it follows $E_4\mct_{1,6k} \subseteq \mct_{1,6k+4}$. Lastly, a) ensures isomorphy between both spaces.
\end{prf}

Theorem \ref{CONG4CASE} is a direct consequence of Proposition \ref{baseprop}.
\begin{prf}[of Theorem \ref{CONG4CASE}]
We combine Proposition \ref{baseprop} b) and the fact that $E_4(z)$ is normalized. Hence, $E_4(z) \ft_{1,6k}(z)$ is a normalized extremal quasimodular form of weight $6k+4$ and depth $1$. Uniqueness of such a function forces that
\begin{align*}
\ft_{1,6k+4}(z) = E_4(z) \ft_{1,6k}(z).
\end{align*}
Recalling $E_4(z) \in \Zb[[q]]$, Theorem \ref{pellthm} implies
\begin{align*}
\ft_{1,6k+4}(z) \in \Zb\left[\frac{1}{p} \colon p < 6k\right][[q]]
\end{align*}
for every $k \geq 0$, which is the statement of Theorem \ref{CONG4CASE}.
\end{prf}

\section{Additional structure in depth $1$} \label{d1sec}
\hspace{\parindent} Henceforth, we omit the dependencies on $z$ and $q$, where it is clear. We recall the normalized modular discriminant function
\begin{align*}
\Delta(z) \coloneqq q\prod_{n=1}^{\infty} \left(1-q^n\right)^{24} = \sum_{n \geq 1} \tau(n)q^n \in \Zb[[q]] \cap \mc_{12}.
\end{align*}
We follow \cite[460-462]{kk06} and define four sequences of polynomials
\begin{align*}
P_0(x) &= 1, \qquad P_1(x) = x, \qquad P_{k+1}(x) = xP_k(x)  + \mu_kP_{k-1}(x), \\
Q_0(x) &= 0, \qquad Q_1(x) = 1, \qquad Q_{k+1}(x) = xQ_k(x)  + \mu_kQ_{k-1}(x), \\
P_0^*(x) &= 1, \qquad P_1^*(x) = x, \qquad P_{k+1}^*(x) = xP_k^*(x)  + \mu^*_kP_{k-1}^*(x), \\
Q_0^*(x) &= 0, \qquad Q^*_1(x) = 1, \qquad Q_{k+1}^*(x) = xQ_k^*(x)  + \mu_k^*Q_{k-1}^*(x), \\
\mu_k &= \frac{12(6k+1)(6k+5)}{k(k+1)}, \qquad \qquad \mu_k^* = \frac{12(6k-1)(6k+7)}{k(k+1)}.
\end{align*}
These polynomials encode an inductive structure of extremal quasimodular forms of depth $1$.

\begin{thm}[\protect{\cite[Theorem 2.1]{kk06}}] \label{kkthm}
\begin{enumerate}[(a)]
\item Suppose that $w=6k$, $k \geq 1$. Then
\begin{align*}
f_{1,6k} \coloneqq \Delta^{\frac{k-1}{2}}P_{k-1}\left(\frac{E_6}{\Delta^{\frac{1}{2}}}\right)\frac{D(E_4)}{240}-\Delta^{\frac{k}{2}}Q_{k-1}\left(\frac{E_6}{\Delta^{\frac{1}{2}}}\right)
\end{align*}
is an extremal quasimodular form of weight $w$ and depth $1$ on $\Gamma$. Furthermore, $f_{1,6k}$ and is a solution of the differential equation
\begin{align*}
D^2(f)-\frac{w}{6}E_2D(f)+\frac{w(w-1)}{12}D(E_2)f = 0.
\end{align*}
\item Suppose that $w=6k+2$, $k \geq 1$. Then
\begin{align*}
f_{1,6k+2} \coloneqq \Delta^{\frac{k-1}{2}}P_{k-1}^*\left(\frac{E_6}{\Delta^{\frac{1}{2}}}\right)\left(-\frac{D(E_6)}{504}\right)-\Delta^{\frac{k}{2}}Q_{k-1}^*\left(\frac{E_6}{\Delta^{\frac{1}{2}}}\right)E_2
\end{align*}
is an extremal quasimodular form of weight $w$ and depth $1$ on $\Gamma$. Furthermore, $f_{1,6k+2}$  is a solution of the differential equation
\begin{align*}
D^2(f)-\left(\frac{w}{6}E_2-\frac{1}{3}\frac{E_6}{E_4}\right)D(f)+\left(\frac{w(w-1)}{12}D(E_2)-\frac{w-1}{18}\frac{D(E_6)}{E_4}\right)f = 0.
\end{align*}
\item An extremal quasimodular form of weight $w \equiv 4 \pmod{6}$ and depth $1$ is obtained from the form in a) with $w$ replaced by $w-4$, by multiplying by $E_4$. The differential equation it satisfies is
\begin{align*}
D^2(f)-\left(\frac{w}{6}E_2-\frac{2}{3}\frac{E_6}{E_4}\right)D(f)+\left(\frac{w(w-1)}{12}D(E_2)-\frac{w-1}{9}\frac{D(E_6)}{E_4}-\frac{2}{9}\left(E_4-\frac{E_6^2}{E_4^2}\right)\right)f = 0.
\end{align*}
\end{enumerate}
\end{thm}

\begin{prf}
The development of this result is contained in the proof of \cite[Theorem 2]{kk03}, be aware of the fact that $w$ is replaced by $w+1$ there. Then inspect dimension and exponents to conclude that the solution is indeed extremal. See \cite[461,462]{kk06} as well.
\end{prf}

Checking the parities of the four polynomials, we observe that the apperances of the square root in Theorem \ref{kkthm} cancel each other.
\begin{cor} \label{kkcor}
For every integer $k \geq 0$ we have
\begin{align*}
f_{1,6(k+2)} &= E_6 f_{1,6(k+1)} + \mu_k \Delta f_{1,6k}, \\
f_{1,6(k+2)+2} &= E_6 f_{1,6(k+1)+2} + \mu_k^* \Delta f_{1,6k+2},
\end{align*}
where we set $\mu_0 \coloneqq \mu_0^* \coloneqq -1$ in accordance with Proposition \ref{baseprop} a).
\end{cor}

\begin{prf}
This is a short argument by induction on $k \geq 0$. The induction step is simply evaluating the polynomial recursions from Theorem \ref{kkthm}.
\end{prf}

On one hand, Example \ref{ex1} and the proof of Theorem \ref{CONG4CASE} demonstrate
\begin{align*}
\ft_{1,2}, \ft_{1,6}, \ft_{1,8}, \ft_{1,10} \in \Zb[[q]].
\end{align*}
In addition, we use the following proposition during the proof of Theorem \ref{INITHM}.
\begin{prop}[G. Nebe, \protect{\cite[11]{pellarin}}]
We have $\ft_{1,14} \in \Zb[[q]]$.
\end{prop}

On the other hand, we invoke the machinery of Corollary \ref{kkcor}.
\begin{exa}
\begin{enumerate}
\item Let $w = 12$. We calculate that
\begin{align*}
f_{1,12} &= E_6f_{1,6} + \mu_0\Delta = E_6\ft_{1,6} - \Delta = -462q^2 -25872q^3  + O\left(q^4\right),
\end{align*}
getting
\begin{align*}
\ft_{1,12}(z) &\in \Zb\left[\frac{1}{p} \colon p < 12\right][[q]],
\end{align*}
since $462 = 2 \cdot 3 \cdot 7 \cdot 11$.
\item Let $w = 20$. We calculate that
\begin{align*}
f_{1,20} &= E_6f_{1,14} + \mu_1^*\Delta f_{1,8} = E_6^2\ft_{1,8} + \mu_1^*\Delta \ft_{1,8} = 163020q^3+29832660q^4+O\left(q^5\right),
\end{align*}
getting
\begin{align*}
\ft_{1,20}(z) &\in \Zb\left[\frac{1}{p} \colon p < 20\right][[q]],
\end{align*}
since $163020 = 2^2 \cdot 3 \cdot 5 \cdot 11 \cdot 13 \cdot 19$.
\end{enumerate}
\end{exa}
However, this does not suffice to prove that the results of Theorem \ref{pellthm} and Theorem \ref{INITHM} are optimal. Towards this direction, we remark that Theorem \ref{kkthm} and Corollary \ref{kkcor} can be translated to a recursive description of $\ft_{1,w}$ too, including a formula for the normalizing coefficient. Both were performed in \cite[462--464]{kk06}. Some numerical experiments suggest that \cite[Theorem 2.3]{kk06} does not suffice to attack both Theorem \ref{INITHM} and its sharpness. Indeed, the formulas provided there do not reveal any divisibility properties possibly shared by all Fourier coefficients.

\section{Proof of Theorem \ref{INITHM}}

\hspace{\parindent} We choose the unique normalized square root in $q^\frac{1}{2}\Cb[[q]]$.

\subsection{Strategy of the proof} \label{stratsubsec}

\hspace{\parindent} The novel observation is that the proof strategy of Pellarin, c.\ f.\ \cite[6-10]{pellarin}, applies to to our setting of weight $w = 6k+2$ and depth $1$ as well. Proposition \ref{odeprop} b) provides an analogue of his starting point \cite[Proposition 3.1]{pellarin}. Then we noted that his remaining intermediate results carry over almost verbatim to our setting, although our starting point b) is more complicated. For the sake of completeness, we state \cite[Proposition 3.1]{pellarin} as item a) of the following Proposition.
\begin{prop} \label{odeprop}
If $f_{1,w}$ denotes one of the extremal quasimodular forms of weight $w$ and depth $1$ from Theorem \ref{kkthm}, then
\begin{enumerate}[(a)]
\item $D^2\left(f_{1,6k}\Delta^{-\frac{k}{2}}\right) = \frac{k^2}{4}E_4\cdot\left(f_{1,6k}\Delta^{-\frac{k}{2}}\right)$,
\item $D^2\left(f_{1,6k+2}\Delta^{-\frac{k}{2}}E_4^{-\frac{1}{2}}\right) = \left(\left(\frac{k^2}{4}-\frac{1}{12}\right)E_4+\frac{1}{12}\left(\frac{E_6}{E_4}\right)^2\right)\cdot\left(f_{1,6k+2}\Delta^{-\frac{k}{2}}E_4^{-\frac{1}{2}}\right)$.
\end{enumerate}
\end{prop}

We will use this result in the following manner.
\begin{cor} \label{odecor}
The forms $\ft_{1,6k}\Delta^{-\frac{k}{2}}$ and $\ft_{1,6k+2}\Delta^{-\frac{k}{2}}E_4^{-\frac{1}{2}}$ are the unique normalized solutions of the corresponding ordinary differential equations given in Proposition \ref{odeprop} a) and b) respectively.
\end{cor}

\begin{prf}[of Corollary \ref{odecor}]
Since all factors are normalized, the forms $\ft_{1,6k}\Delta^{-\frac{k}{2}}$ and $\ft_{1,6k+2}\Delta^{-\frac{k}{2}}E_4^{-\frac{1}{2}}$ are normalized. By uniqueness, which normalization guarantees, $\ft_{1,6k}$ and $\ft_{1,6k+2}$ coincide with the normalized forms of $f_{1,6k}$ and $f_{1,6k+2}$ from Theorem \ref{kkthm}. Finally, normalization does not affect the property of solvability of an ordinary differential equation.
\end{prf}

\begin{rmk} \label{odermk}
\begin{enumerate}
\item Propsition \ref{odeprop} b) may be rewritten as follows. The representation $\Delta = \frac{E_4^3-E_6^2}{1728}$ yields 
\begin{align*}
-\frac{1}{12}\left(E_4-\frac{E_6^2}{E_4^2}\right) = -144\frac{\Delta}{E_4^2}.
\end{align*}
\item The presence of the additional factor $E_4^{-\frac{1}{2}}$ in b) is justified by inspection of weights. It compensates for the change from $w = 6k$ to $w = 6k+2$. For this purpose, we can not take any nontrivial power of $E_2$ without manipulating the depth, which is fixed to $1$.  Hence, there are just two candidates to consider: $E_4^{-\frac{1}{2}}$ and $E_6^{-\frac{1}{3}}$. Recalling the result $E_4\mct_{6k}^1 = \mct_{6k+4}^1$ justifies the investigation of the contribution of $E_4^{-\frac{1}{2}}$ instead of $E_6^{-\frac{1}{3}}$. Consequently, the choice of $E_4^{-\frac{1}{2}}$ appears to be natural in our setting. Summing up, we preserve an ordinary differential equation of the form
\begin{align*}
D^2(g) = h(E_4, E_6)g,
\end{align*}
where $g$ has weight $0$ in both cases and the term $h(E_4, E_6)$ has to have weight $4$, since $D^2(g)$ has weight $4$. Once more, note that $h$ has to be independent of $E_2$ and that we may omit a separate dependency of $h$ on $\Delta$.
\end{enumerate}
\end{rmk}

Following Pellarin's approach, define
\begin{align*}
\psi(k) &\coloneqq k+1, \\
\Psi(\mu(k)) &\coloneqq \psi^2 + \mu(k)\left(\frac{E_6}{\Delta^{\frac{1}{2}}}\psi - 1\right).
\end{align*}
Furthermore, we abbreviate
\begin{align*}
\varphi \coloneqq \ft_{1,6k+2}\Delta^{-\frac{k}{2}}E_4^{-\frac{1}{2}}.
\end{align*}
We omit the dependency of $\varphi$ on $k$ because we will not require it.

We investigate the action of the operator $\Psi$ on $\varphi$. Roughly speaking, its purpose is to isolate $\mu_k^*$ from the recursion in Corollary \ref{kkcor}, due to the following Proposition.
\begin{prop} \label{lemma3.7}
We have that $\Psi((\mu_{k+1}^*)^{-1})(\varphi) = 0$.
\end{prop}

The proof is adapted from \cite[9]{pellarin} and relies on three technical lemmas, which are also stated in \cite[8-9]{pellarin}. We postpone their purely computational proofs to Subsection \ref{compsubsec}.

\begin{lemma} \label{lemma3.6}
Recall $\mu_k^* = \frac{12(6k-1)(6k+7)}{k(k+1)}$ appearing in the definition of $P_{k+1}^*$ and $Q_{k+1}^*$. Write formally
\begin{align*}
\Psi(\mu(k))\left(\varphi\right)(z) = q^{\frac{k}{2}}\sum_{n \geq 0} \beta_n(k)q^n.
\end{align*}
Then $\beta_0(k) = 0$. Moreover, $\beta_1(k) = 0$ if and only if $\mu(k) = (\mu_{k+1}^*)^{-1}$.
\end{lemma}

\begin{lemma} \label{lemma3.5}
Define
\begin{align*}
\Dc &\coloneqq D^2 - \frac{k^2}{4}E_4 + 144\frac{\Delta}{E_4^2}, \\
\Fc &\coloneqq \Delta^{-\frac{1}{2}}\big(12E_4D+E_2E_4+(6k+5)E_6\big)\psi.
\end{align*}
Then we have for every choice of $\mu(k)$
\begin{align*}
\psi^2\Dc\psi^{-2}\Psi(\mu(k)) - \Psi(\mu(k))\Dc = \mu(k)E_4\left((k+1)-\frac{\Fc}{12}\right).
\end{align*}
\end{lemma}
Note that, according to Remark \ref{odermk}, $\Dc$ is just the differential operator encoding the result of Proposition \ref{odeprop} b). That is, we have $\Dc\varphi = 0$.

\begin{lemma} \label{lemma3.4}
The operators $\Dc$ and $\Fc$, introduced in Lemma \ref{lemma3.5}, satisfy the commutation rule
\begin{align*}
\Dc\Fc - \Fc\Dc = -4\Delta^{-\frac{1}{2}}\left(E_2E_4+2E_6\right)\psi\Dc.
\end{align*}
\end{lemma}

\begin{rmk} \label{correctionremark}
Our Lemmas \ref{lemma3.5} and \ref{lemma3.4} are the same statements as Pellarin's Lemmas 3.5 and 3.4 except the definition of the operator $\Fc$. Pellarin uses the operator 
\begin{align*}
F \coloneqq \Delta^{-\frac{1}{2}}\left(12E_2D+E_2E_4+(6k+5)E_6\right)\psi
\end{align*}
instead and proposes the same formulas with our $\Fc$ replaced by his $F$. However, our computations proving both of our lemmas are independent of the residue class of the weight. Hence, we suggest to correct his operator $F$ to our operator $\Fc$ in \cite{pellarin} too.
\end{rmk}

Assuming these Lemmas, we are able to prove Proposition \ref{lemma3.7}.
\begin{prf}[of Proposition \ref{lemma3.7}]
We define 
\begin{align*}
V \coloneqq \text{Ker}(\Dc) \cap Y\Qb(k)\left(\left(q^{\frac{1}{2}}\right)\right),
\end{align*}
where $Y$ is a formal solution of $D(Y) = \frac{k}{2}Y$. We proceed as follows. \\
\textbf{Claim 1:} The vector space $V$ is one-dimensional and generated by $\varphi$ over $\Qb(k)$.
\begin{prf}[of Claim 1]
Indeed, 
\begin{align*}
\varphi(z) = Y\sum_{n\geq0} \alpha_n(k)q^n \in V
\end{align*}
which is confirmed in Subsection \ref{compsubsec}. Additionally, the second linearly independent solution to Proposition \ref{odeprop} b) is given by 
\begin{align*}
Y^{-1}\sum_{n\geq0} \alpha_n(-k)q^n
\end{align*}
and thus not contained in $V$. Otherwise, $Y$ and $Y^{-1}$ would be linear dependent over $\Qb(k)\big((q^{\frac{1}{2}})\big)$. Since $V$ is at most two-dimensional, we obtain
\begin{align*}
V = \Qb(k)\varphi
\end{align*}
proving the first claim.
\end{prf}

Utilizing Lemma \ref{lemma3.4} and Proposition \ref{odeprop} b), we deduce that $\Fc(\varphi) \in \text{Ker}(\Dc)$. Thus, by definition of $\Fc$, we get $\Fc(\varphi) \in V$ and we write $\Fc(\varphi) = \lambda \varphi$ for some $\lambda \in \Qb(k)$ according to claim 1. \\
\textbf{Claim 2:} We have $\lambda = 12(k+1)$.
\begin{prf}[of Claim 2]
This is a computation and thus moved to Section \ref{compsubsec} as well.
\end{prf}
Combining Lemma \ref{lemma3.5} and the second claim shows
\begin{align*}
\psi^2\Dc\psi^{-2}\Psi(\mu(k))(\varphi) = 0,
\end{align*}
that is
\begin{align*}
\psi^{-2}\Psi(\mu(k))(\varphi) \in \text{Ker}(\Dc).
\end{align*}
By definition of $\Psi(\mu(k))$, it follows that
\begin{align*}
\psi^{-2}\Psi(\mu(k))(\varphi) \in V,
\end{align*}
which yields 
\begin{align*}
\Psi(\mu(k))(\varphi) \in \Qb(k)\psi^2(\varphi)
\end{align*}
in view of the first claim. Summing up, we obtain
\begin{align*}
\Psi(\mu(k))(\varphi) = \xi\psi^2(\varphi)
\end{align*}
for some $\xi \in \Qb(k)$.
Choosing $\mu(k) = (\mu_{k+1}^*)^{-1}$ and invoking Lemma \ref{lemma3.6}, we deduce that $\xi = 0$. This proves the proposition.
\end{prf}

A similar result holds in the case of weight $w = 6k$ (Pellarin's case):
\begin{align*}
\Psi\left(\left(\mu_{k+1}\right)^{-1}\right)\left(\ft_{1,6k}\Delta^{-\frac{k}{2}}\right) = 0.
\end{align*}
Regarding its proof, one has to replace $\mu_{k+1}^*$ by $\mu_{k+1}$ in Lemma \ref{lemma3.6}, $\Dc$ by $D^2-\frac{k^2}{4}$ in both Lemma \ref{lemma3.5} and Lemma \ref{lemma3.4}, and $\varphi$ by $\ft_{1,6k}\Delta^{-\frac{k}{2}}$. The $q$-expansion of $\ft_{1,6k}\Delta^{-\frac{k}{2}}$ is of the same structure as the $q$-expansion of $\varphi$, since Proposition \ref{odeprop} a) and b) share the term $\frac{k^2}{4}E_4$. We conclude that we may adapt Pellarin's final proof of Theorem \ref{pellthm} to prove our Theorem \ref{INITHM} as well.

\begin{prf}[of Theorem \ref{INITHM}]
The result follows by induction on $j \geq 1$. We remarked in Section \ref{d1sec} that $\ft_{1,8}$, $\ft_{1,14} \in \Zb[[q]]$, so Theorem \ref{INITHM} obviously holds for these two functions. We assume $j \geq 3$. Rewriting Proposition \ref{lemma3.7} yields
\begin{align*}
\varphi - \frac{E_6}{\Delta^{\frac{1}{2}}}\psi(\varphi)  = \mu(k)^{-1}\psi^2(\varphi) = \mu_{k+1}^*\psi^2(\varphi).
\end{align*}
We apply $\psi^{j-2}$ to this equation and evaluate at $k=0$. Formally, writing
\begin{align*}
\left(\psi^{j-2}(\varphi) - \frac{E_6}{\Delta^{\frac{1}{2}}}\psi^{j-1}(\varphi)\right)\Big\vert_{k=0}(z) &\eqqcolon q^{\frac{j}{2}}\sum_{n \geq 0} \iota_n(j) q^n, \\
\left(\psi^j(\varphi)\right)\Big\vert_{k=0}(z) &\eqqcolon q^{\frac{j}{2}}\sum_{n \geq 0} \kappa_n(j) q^n,
\end{align*}
we obtain
\begin{align*}
\kappa_n(j) = (\mu_{j-1}^*)^{-1}\iota_n(j) = \frac{j(j-1)}{12(6j-7)(6j+1)}\iota_n(j).
\end{align*}
The induction hypothesis is the claim of Theorem \ref{INITHM} for all $\psi^m\varphi\vert_{k=0}$, where $m < j$. Equivalently, we suppose that the primes dividing the denominators of $\iota_n(m)$ are smaller than $6m+2$ for every $n$ and $m < j$. Then it is apparent that the additional prime factors caused by the denominator of $ (\mu_{j-1}^*)^{-1}$ have to be smaller than $6j+2$. We conclude that if $p$ is a prime dividing the denominator of $\kappa_n(j)$ then $p < 6j+2$. Finally, 
\begin{align*}
\left(\psi^j(\varphi)\right)\Big\vert_{k=0} = \ft_{1,6j+2}\Delta^{-\frac{j}{2}}E_4^{-\frac{1}{2}}.
\end{align*}
We checked $\Delta^{\frac{1}{2}} \in \Zb[[q]]$ and $E_4^{\frac{1}{2}} \in \Zb[[q]]$ in Lemma \ref{sqrtlemma}. Thus, Theorem \ref{INITHM} holds for $\ft_{1, 6j+2}$. This completes the proof.
\end{prf}

\begin{rmk}
For the proof of Theorem \ref{pellthm}, one may simply replace $\varphi$ by $\ft_{1,6k}\Delta^{-\frac{k}{2}}$ and $\mu^*$ by $\mu$. Then, observe that
\begin{align*}
\mu_{j-1}^{-1} = \frac{j(j-1)}{12(6j-5)(6j-1)},
\end{align*}
and hence it suffices to rephrase the induction hypthesis suitably to conclude in a similar way.
\end{rmk}

\subsection{Computations} \label{compsubsec}
\hspace{\parindent} We collect all loose ends, namely we provide the missing, entirely computational proofs. We begin with ``completing'' Ramanujan's differential system by the following fruitful result.
\begin{lemma} \label{deltadiff}
We have $D(\Delta) = E_2\Delta$.
\end{lemma}

\begin{prf}[of Lemma \ref{deltadiff}]
This is a straightforward observation using $\Delta = \frac{E_4^3-E_6^2}{1728}$ and Ramanujan's differential system.
\end{prf}

This enables us to verify Proposition \ref{odeprop}.
\begin{prf}[of Proposition \ref{odeprop}]
Two proofs of a) were given in \cite[Propsition 3.1]{pellarin}. A third proof is purely computational and relies on Lemma \ref{deltadiff} and Ramanujan's differential system. First, we compute
\begin{align}
D^2(fg^{\alpha}) = g^{\alpha}\left[D^2(f) + 2\alpha \frac{D(g)}{g}D(f) + \left(\alpha(\alpha-1)\left(\frac{D(g)}{g}\right)^2+\alpha\frac{D^2(g)}{g}\right)f\right].
\end{align}
Inserting $\alpha=\frac{1}{2}$, $g=\Delta^{-k}$, $E_2^2 = 12D(E_2) + E_4$ into $(1)$ and collecting all terms gives directly
\begin{align*}
D^2\left(f_{1,6k}\Delta^{-\frac{k}{2}}\right) = \Delta^{-\frac{k}{2}}\left(D^2(f_{1,6k})-\frac{w}{6}E_2D(f_{1,6k})+\frac{w(w-1)}{12}D(E_2)f_{1,6k}\right) + \frac{k^2}{4}E_4\cdot\left(f_{1,6k}\Delta^{-\frac{k}{2}}\right)
\end{align*}
and item a) follows according to theorem \ref{kkthm} a).

The calculation verifying b) is a little bit more tedious and we provide some intermediate steps. Inserting $\alpha=\frac{1}{2}$, $g=\Delta^{-k}E_4^{-1}$, $D(E_4) = \frac{E_2E_4-E_6}{3}$ into $(1)$ and collecting all terms gives
\begin{align*} 
D^2(f\Delta^{-\frac{k}{2}}E_4^{-\frac{1}{2}}) = \Delta^{-\frac{k}{2}}E_4^{-\frac{1}{2}} \cdot & \Bigg\{D^2(f) + \left(-\frac{w}{6}E_2+\frac{1}{3}\frac{E_6}{E_4}\right)D(f) \\
& + \left[\frac{1}{4}\left(-\frac{w}{6}E_2+\frac{1}{3}\frac{E_6}{E_4}\right)^2+\frac{1}{2}\left(-\frac{w}{6}D(E_2)+\frac{1}{3}D\left(\frac{E_6}{E_4}\right)\right)\right]f\Bigg\}.
\end{align*}
Utilizing Ramanujan's differential system, we simplify
\begin{align*}
& \frac{1}{4}\left(-\frac{w}{6}E_2+\frac{1}{3}\frac{E_6}{E_4}\right)^2+\frac{1}{2}\left(-\frac{w}{6}D(E_2)+\frac{1}{3}D\left(\frac{E_6}{E_4}\right)\right) \\
&= \frac{w^2}{144}(12D(E_2)+E_4) + \left(\frac{1}{36}+\frac{1}{18}\right)\frac{E_6^2}{E_4^2} -\frac{w}{12}D(E_2) + \left(-\frac{w}{36}-\frac{1}{18}\right)\frac{E_2E_6-E_4^2+E_4^2}{E_4} + \frac{1}{6}\frac{D(E_6)}{E_4} \\
&= \frac{w(w-1)}{12}D(E_2) + \frac{-w+1}{18}\frac{D(E_6)}{E_4} + \left(\left(\frac{w-2}{12}\right)^2-\frac{1}{12}\right)E_4 + \frac{1}{12}\frac{E_6^2}{E_4^2}
\end{align*}
and the result follows by virtue of Theorem \ref{kkthm} b).
\end{prf}

Next, we compute the $q$-expansions of the functions involved in Proposition \ref{odeprop} and confirm $\varphi \in V$, where $V$ is defined in the proof of Proposition \ref{lemma3.7}. First, let
\begin{align*}
\left(\frac{k^2}{4}-\frac{1}{12}\right)E_4(z)+\frac{1}{12}\left(\frac{E_6}{E_4}\right)^2(z) \eqqcolon \sum_{n \geq 0} a_n(k)q^n
\end{align*}
Then
\begin{align*}
\frac{1}{E_4}(z) &= 1 + \sum_{n \geq 1} b_nq^n \ , \ b_0 = 1 \ , \ b_n = -\sum_{m=1}^n 240\sigma_3(m)b_{n-m}, \\
\frac{E_6}{E_4}(z) &= 1 + \sum_{n \geq 1} c_nq^n \ , \ c_0 = 1 \ , \ c_n = -504\sigma_5(n)-\sum_{m=1}^n 240\sigma_3(m)c_{n-m}, \\
\left(\frac{E_6}{E_4}\right)^2(z) &= 1 + \sum_{n \geq 1} d_nq^n \ , \ d_0 = 1 \ , \ d_n = 2c_n + \sum_{m=1}^n c_mc_{n-m}, \\
a_0(k) &= \frac{k^2}{4} \in \Zb \big[\frac{1}{2}\big] \ , \ \forall n\geq 1 \colon a_n(k) = 240\sigma_3(n)\left(\frac{k^2}{4}-\frac{1}{12}\right) + \frac{d_n}{12} \in \Zb.
\end{align*}
Second, let
\begin{align*}
\varphi(z) = q^{\frac{k}{2}}\sum_{n\geq0} \alpha_n(k)q^n
\end{align*}
where $q^{\frac{k}{2}}$ will take care of $a_0(k) = \frac{k^2}{4}$. We insert this $q$-expansion into Proposition \ref{odeprop} b) giving
\begin{align*}
\alpha_n(k) = \frac{1}{n(n+k)}\sum_{m=1}^n a_m(k)\alpha_{n-m}(k)
\end{align*}
for every $n \geq 1$. We may choose $\alpha_0(k) = 1$ for every $k$.

\begin{lemma} \label{sqrtlemma}
We have that $\Delta^{\frac{1}{2}} \in \Zb[[q]]$ and that $E_4^{\frac{1}{2}} \in \Zb[[q]]$.
\end{lemma}

\begin{prf}
Writing formally
\begin{align*}
\sum_{n \geq N} r_n q^n = \left(\sum_{n \geq N} s_n q^n\right)^2
\end{align*}
yields
\begin{align*}
s_N = \sqrt{r_N} \qquad s_n = \frac{1}{s_N}\left(r_n - \sum_{m=N+1}^n s_m s_{n-m}\right)
\end{align*}
for all $n\geq N+1$. Applying this formula to $E_4 \in \Zb[[q]]$ with $N=0$, $r_0 = 1$, $r_n = 240\sigma_3(n) $ and to $\Delta \in \Zb[[q]]$ with $N=1$, $r_1 = 1$, $r_n = \tau(n) $, we deduce the claim by recursion.
\end{prf}

The previous Lemma completes the last step in the proof of Theorem \ref{INITHM}. Furthermore, we use it in the following special case.
\begin{lemma} \label{deltanegroot}
The $q$-expansion of $\Delta^{-\frac{1}{2}}$ begins as follows: 
\begin{align*}
\Delta^{-\frac{1}{2}}(z) = q^{-\frac{1}{2}} + 12q^{\frac{1}{2}} + O\left(\left(q^{\frac{1}{2}}\right)^3\right)
\end{align*}
\end{lemma}

\begin{prf}[of Lemma \ref{deltanegroot}]
This follows by the formula for the square root of a $q$-expansion from the previous Lemma and by the inversion formula of a $q$-expansion analogously to $\frac{1}{E_4}$ above. Alternatively, one may argue using the representation of $\Delta$ as an infinite product.
\end{prf}

We collect everything together to prove Lemma \ref{lemma3.6}.
\begin{prf}[of Lemma \ref{lemma3.6}]
Clearly we have
\begin{align*}
\psi^2\left(q^{\frac{k}{2}}\sum_{n\geq0} \alpha_n(k)q^n\right) &= q^{\frac{k+2}{2}}\sum_{n\geq0} \alpha_n(k+2)q^n, \\
\mu(k)\left(\frac{E_6}{\Delta^{\frac{1}{2}}}\psi - 1\right)\left(q^{\frac{k}{2}}\sum_{n\geq0} \alpha_n(k)q^n\right) &= \mu(k)E_6\Delta^{-\frac{1}{2}}q^{\frac{k+1}{2}}\sum_{n\geq0} \alpha_n(k+1)q^n - \mu(k)q^{\frac{k}{2}}\sum_{n\geq0} \alpha_n(k)q^n.
\end{align*}
Hence, according to Lemma \ref{deltanegroot},
\begin{align*}
\beta_0(k) &= \mu(k) - \mu(k) = 0, \\
\beta_1(k) &= \alpha_0(k+2) + \mu(k)\left(\alpha_1(k+1)+(12-504)\alpha_0(k+1)-\alpha_1(k)\right) \\
&= 1 + \mu(k)\left(\frac{60(k+1)^2-144}{k+2} - 492 - \frac{60k^2-144}{k+1}\right) \\
&= 1-\mu(k)\mu_{k+1}^*
\end{align*}
and the result follows.
\end{prf}

We now turn our interest towards operator calculus and state the following rules to distinguish between (non-commutative) multiplication by an operator and application of an operator. Let $g \in \Qb(k)\big((q^{\frac{1}{2}})\big)$. Then, we have
\begin{align*}
D\psi = \psi D, \qquad Dg = gD + D(g), \qquad \psi g = \psi(g)\psi.
\end{align*}

With this established, recall 
\begin{align*}
\Dc &= D^2 - \frac{k^2}{4}E_4 + 144\frac{\Delta}{E_4^2}, \\
\Fc &= \Delta^{-\frac{1}{2}}\big(12E_4D+E_2E_4+(6k+5)E_6\big)\psi.
\end{align*}

We are now in position to prove the Lemmas \ref{lemma3.5} and \ref{lemma3.4}.
\begin{prf}[of Lemma \ref{lemma3.5}]
Substituting
\begin{align*}
\Dc \eqqcolon D^2-\rho, \quad \Psi(\mu(k)) \eqqcolon \psi^2 + \eta(k),
\end{align*}
we calculate
\begin{align*}
\psi^2\Dc\psi^{-2}\Psi(\mu(k)) - \Psi(\mu(k))\Dc = D^2\eta(k)-\eta(k)D^2 - \psi^2\rho\psi^{-2}\eta(k) + \eta(k)\rho.
\end{align*}
Substituting back, we have
\begin{align*}
& \psi^2\Dc\psi^{-2}\Psi(\mu(k)) - \Psi(\mu(k))\Dc \\
& = \mu(k)\left\{\left(D^2\frac{E_6}{\Delta^{\frac{1}{2}}}\psi-\frac{E_6}{\Delta^{\frac{1}{2}}}\psi D^2\right) + \left(\frac{(k+1)^2}{4}\frac{E_6}{\Delta^{\frac{1}{2}}}E_4\psi - \frac{(k+2)^2}{4}\frac{E_6}{\Delta^{\frac{1}{2}}}E_4\psi\right) + \left(\frac{(k+2)^2}{4}-\frac{k^2}{4}\right)E_4\right\},
\end{align*}
since 
\begin{align*}
144\frac{\Delta}{E_4^2}\eta(k) - 144\frac{\Delta}{E_4^2}\eta(k) = 0.
\end{align*}
Now, by the aforementioned rules
\begin{align*}
D^2g = DDg = DgD+DD(g) = gD^2 + 2D(g)D + D^2(g).
\end{align*}
Moreover, utilizing Ramanujan's differential system and $D(\Delta) = E_2\Delta$, we have that
\begin{align*}
D\left(\Delta^{-\frac{1}{2}}E_6\right) &= -\frac{1}{2}\Delta^{-\frac{1}{2}} E_4^2, \\
D^2\left(\Delta^{-\frac{1}{2}}E_6\right)&= \frac{1}{3}\Delta^{-\frac{1}{2}}E_4E_6-\frac{1}{12}\Delta^{-\frac{1}{2}}E_2E_4^2.
\end{align*}
Inserting $g = \Delta^{-\frac{1}{2}}E_6$ and combining, this produces
\begin{align*}
& \psi^2\Dc\psi^{-2}\Psi(\mu(k)) - \Psi(\mu(k))\Dc \\
& = \mu(k)\left\{-\Delta^{-\frac{1}{2}}E_4^2D\psi+\left(\frac{1}{3}-\frac{(k+2)^2}{4}+\frac{(k+1)^2}{4}\right)\Delta^{-\frac{1}{2}}E_4E_6-\frac{1}{12}\Delta^{-\frac{1}{2}}E_2E_4^2\right\}\psi + \mu(k)(k+1)E_4,
\end{align*}
and we are done up to rearrangement of the terms.
\end{prf}

\begin{prf}[of Lemma \ref{lemma3.4}]
We correct and restate Pellarin's Lemma 3.4 in virtue of Remark \ref{correctionremark} as follows.
\begin{align*}
\left(D^2-\frac{k^2}{4}E_4\right)\Fc - \Fc\left(D^2-\frac{k^2}{4}E_4\right) = -4\Delta^{-\frac{1}{2}}\left(E_2E_4+2E_6\right)\psi\left(D^2-\frac{k^2}{4}E_4\right)
\end{align*}
\textbf{Step 1:} We prove this equation first: \\
Set $g \coloneqq \Delta^{-\frac{1}{2}}\left(E_2E_4+(6k+5)E_6\right)$. Then, using $D^2g = gD^2 + 2D(g)D + D^2(g)$, we see that
\begin{align*}
D^2\Fc - \Fc D^2 = \left(24D\left(\Delta^{-\frac{1}{2}}E_4\right)D^2+12D^2(\Delta^{-\frac{1}{2}}E_4)D+2D(g)D+D^2(g)\right)\psi.
\end{align*}
Observe that we obtain the first term on the right hand side, namely
\begin{align*}
24D\left(\Delta^{-\frac{1}{2}}E_4\right)D^2\psi = -4\Delta^{-\frac{1}{2}}\left(E_2E_4+2E_6\right)\psi D^2
\end{align*}
The second term requires more effort. On one hand,
\begin{align*}
&-\frac{k^2}{4}E_4\Fc + \Fc\frac{k^2}{4}E_4 = -3k^2\Delta^{-\frac{1}{2}}E_4^2D\psi + 3(k+1)^2\Delta^{-\frac{1}{2}}E_4(E_4D + D(E_4))\psi + \frac{2k+1}{4}E_4g\psi \\
&= (k+1)^2\Delta^{-\frac{1}{2}}E_4\left(E_2E_4-E_6\right)\psi + (2k+1)\Delta^{-\frac{1}{2}}\left(3E_4^2D+\frac{1}{4}\Big(E_2E_4^2+(6k+5)E_4E_6\Big)\right)\psi.
\end{align*}
On the other hand, focussing on the $k$-dependent term of $D^2(g)$ for the moment,
\begin{align*}
D^2\left((6k+5)\Delta^{-\frac{1}{2}}E_6\right) = (6k+5)\Delta^{-\frac{1}{2}}\left(-\frac{1}{12}E_2E_4^2+\frac{1}{3}E_4E_6\right).
\end{align*}
Summing up the previous two equations, we obtain
\begin{align*}
& -\frac{k^2}{4}E_4\Fc + \Fc\frac{k^2}{4}E_4 + D^2\left((6k+5)\Delta^{-\frac{1}{2}}E_6\right)\psi \\
&= 3(2k+1)\Delta^{-\frac{1}{2}}E_4^2D\psi + \Delta^{-\frac{1}{2}}\left((k+1)^2E_4\left(E_2E_4-E_6\right) + 3(k+1)^2E_4E_6 -\frac{1}{12}E_4E_6 - \frac{1}{6}E_2E_4^2\right)\psi \\
&= 4\Delta^{-\frac{1}{2}}\left(E_2E_4+2E_6\right)\psi\frac{k^2}{4}E_4 + \Delta^{-\frac{1}{2}}\left(3(2k+1)E_4^2D -\frac{1}{12}E_4E_6 - \frac{1}{6}E_2E_4^2\right)\psi
\end{align*}
which establishes the second term on the right hand side of the first step. Thus, collecting the remaining terms, we seek to show that
\begin{align*}
12D^2\left(\Delta^{-\frac{1}{2}}E_4\right)D + 2D(g)D + D^2\left(\Delta^{-\frac{1}{2}}E_2E_4\right) + \Delta^{-\frac{1}{2}}\left(3(2k+1)E_4^2D -\frac{1}{12}E_4E_6 - \frac{1}{6}E_2E_4^2\right) = 0.
\end{align*}
To check this claim, we compute
\begin{align*}
D^2(\Delta^{-\frac{1}{2}}E_4) &= \Delta^{-\frac{1}{2}}\left(\frac{1}{12}E_2^2E_4-\frac{1}{72}E_4(E_2^2-E_4)-\frac{1}{18}E_2(E_2E_4-E_6)+\frac{1}{6}E_4^2\right) \\
12D^2(\Delta^{-\frac{1}{2}}E_4)D &= \Delta^{-\frac{1}{2}}\left(\frac{1}{6}E_2^2E_4+\frac{13}{6}E_4^2+\frac{2}{3}E_2E_6\right)D \\
D(g) &= \Delta^{-\frac{1}{2}}\Big(-\frac{1}{2}E_2(E_2E_4+(6k+5)E_6)+\frac{1}{12}E_4(E_2^2-E_4) \\
& \qquad \qquad +\frac{1}{3}E_2(E_2E_4-E_6)+\frac{6k+5}{2}(E_2E_6-E_4^2)\Big) \\
2D(g)D &= \Delta^{-\frac{1}{2}}\left(-\frac{1}{6}E_2^2E_4+\left(-\frac{1}{6}-(6k+5)\right)E_4^2-\frac{2}{3}E_2E_6\right)D
\end{align*}
which verifies that
\begin{align*}
12D^2\left(\Delta^{-\frac{1}{2}}E_4\right)D + 2D(g)D + 3(2k+1)\Delta^{-\frac{1}{2}}E_4^2D = 0.
\end{align*}
Lastly,
\begin{align*}
D^2\left(\Delta^{-\frac{1}{2}}E_2E_4\right) &= \Delta^{-\frac{1}{2}}\Bigg\{\frac{1}{24}\Big(E_2^3E_4+E_2E_4^2\Big)+\frac{1}{6}E_2^2E_6-\frac{1}{72}E_2E_4\Big(E_2^2-E_4\Big) -\frac{1}{36}E_2^2(E_2E_4-E_6) \\
& \qquad \qquad -\frac{1}{18}E_4(E_2E_4-E_6)-\frac{1}{36}E_6\Big(E_2^2-E_4\Big)-\frac{1}{6}E_2\Big(E_2E_6-E_4^2\Big) \Bigg\} \\
&= \Delta^{-\frac{1}{2}}\left(\frac{1}{6}E_2E_4^2+\frac{1}{12}E_4E_6\right),
\end{align*}
and thus
\begin{align*}
D^2\left(\Delta^{-\frac{1}{2}}E_2E_4\right) + \Delta^{-\frac{1}{2}}\left(-\frac{1}{12}E_4E_6 - \frac{1}{6}E_2E_4^2\right) = 0
\end{align*}
as desired. This proves step $1$.

\textbf{Step 2:} Now, linearity of Step $1$ and $\Dc = D^2-\frac{k^2}{4}E_4 + 144\frac{\Delta}{E_4^2}$ yields
\begin{align*}
\Dc\Fc - \Fc\Dc = -4\Delta^{-\frac{1}{2}}\left(E_2E_4+2E_6\right)\psi\left(D^2-\frac{k^2}{4}E_4\right) + 144\left(\frac{\Delta}{E_4^2}\Fc - \Fc\frac{\Delta}{E_4^2}\right).
\end{align*}
We use the definition of $\Fc$, cancel the commuting terms and then the fact that $Dh = hD+D(h)$ to obtain
\begin{align*}
\frac{\Delta}{E_4^2}\Fc - \Fc\frac{\Delta}{E_4^2} &= 12\Delta^{-\frac{1}{2}}E_4\left(\frac{\Delta}{E_4^2}D - D\frac{\Delta}{E_4^2}\right)\psi = -12\Delta^{-\frac{1}{2}}E_4D\left(\frac{\Delta}{E_4^2}\right)\psi \\
&=  -12\Delta^{-\frac{1}{2}}E_4\frac{\Delta E_2E_4^2-2\Delta E_4\frac{E_2E_4-E_6}{3}}{E_4^4}\psi = -4\Delta^{-\frac{1}{2}}\left(E_2E_4+2E_6\right)\frac{\Delta}{E_4^2}\psi \\
&= -4\Delta^{-\frac{1}{2}}\left(E_2E_4+2E_6\right)\psi\frac{\Delta}{E_4^2},
\end{align*}
which verifies the claim of Lemma \ref{lemma3.4}.
\end{prf}

To finish, we provide the missing calculation postponed during the proof of Proposition \ref{lemma3.7}.
\begin{prf}[of Proposition \ref{lemma3.7}, Claim 2]
We established $\Fc(\varphi) = \lambda \varphi$ for some $\lambda \in \Qb(k)$ during the proof of Proposition \ref{lemma3.7}. We now seek to show that $\lambda = 12(k+1)$. Let
\begin{align*}
\Fc(\varphi)(z) \eqqcolon q^{\frac{k}{2}}\sum_{n\geq 0} \omega_n(k)q^n
\end{align*}
Hence, $\omega_n(k) = \lambda\alpha_n(k)$ and the claim reduces to $\omega_0(k) = 12(k+1)$ due to the fact that $\alpha_0(k) = 1$ for every $k$. We have
\begin{align*}
\Fc(\varphi) &= \Delta^{-\frac{1}{2}}\big(12E_4D+E_2E_4+(6k+5)E_6\big)\psi(\varphi), \\
\psi(\varphi) &=  q^{\frac{k+1}{2}}\sum_{n\geq 0} \alpha_n(k+1)q^n, \\
D\psi(\varphi) &= q^{\frac{k}{2}}\left(\frac{k+1}{2}q^{\frac{1}{2}}\sum_{n\geq 0} \alpha_n(k+1)q^n + q^{\frac{1}{2}}\sum_{n\geq 1} n\alpha_n(k+1)q^n\right).
\end{align*}
In addition, according to Lemma \ref{deltanegroot}, we recall that
\begin{align*}
\Delta^{-\frac{1}{2}}(z) &= q^{-\frac{1}{2}} + 12q^{\frac{1}{2}} + O\left(\left(q^{\frac{1}{2}}\right)^3\right),
\end{align*}
and that that the constant coefficient of $E_2$, $E_4$, and $E_6$ is equal to $1$. Thus, collecting all coefficients contributing to the term $q^{\frac{k}{2}}$ in $\Fc(\varphi)$, we obtain
\begin{align*}
\omega_0(k) = \left(12\frac{k+1}{2} + 1 + (6k+5)\right) = 12(k+1)
\end{align*}
as claimed.
\end{prf}

This result completes the proof of Theorem \ref{INITHM}.


\section{Outlook} \label{outsec}
\hspace{\parindent} In \cite[465-468]{kk06}, Kaneko and Koike described an analogous result to Theorem \ref{kkthm} in the case of weight divisible by $4$ and in depth $2$. In other words, there is an inductive structure of certain extremal quasimodular forms of weight divisible by $4$ and of depth $2$, which is of a similar nature as described in Theorem \ref{kkthm} a). In addition, those forms solve a particular ordinary differential equation again. Both observations together constitued the basis of the proofs of both Theorem \ref{pellthm} and Theorem \ref{INITHM}. Therefore, a suitable rephrasement of 
\cite[6-10]{pellarin} and our Sections \ref{d1sec}, \ref{stratsubsec} might provide an alternative proof of Grabner's result
\begin{align*}
\ft_{2,4k}(z) \in \Zb\left[\frac{1}{p} \colon p < 4k\right][[q]]
\end{align*}
in this special case.

\bigskip

\printbibliography

\bigskip

\bigskip

\begin{flushright}
\textsc{Andreas Mono \\
University of Cologne \\
Department of Mathematics and Computer Science \\
Division of Mathematics \\
Weyertal 86-90 \\
50931 Cologne, Germany} \\
amono@math.uni-koeln.de
\end{flushright}


\end{document}